\documentclass[12pt]{article}
\usepackage{setspace}
\usepackage{amsfonts}
\usepackage{amsmath}
\usepackage{amssymb}
\usepackage{amsthm}
\usepackage{mathrsfs}
\usepackage[all]{xy}

\newcommand{\proofend}{\hfill \hbox{\vrule width 5pt height 5pt depth
0pt}}
\newcommand{\R}{\mathbb{R}}
\newcommand{\C}{\mathbb{C}}
\newcommand{\Z}{\mathbb{Z}}
\newcommand{\N}{\mathbb{N}}
\newcommand{\Q}{\mathbb{Q}}

\newcommand{\Oz}{\mathcal{O}}
\newcommand{\ord}{\mathrm{ord}}

\newcommand{\proj}{\mathbb{P}}
\newcommand{\A}{\mathbb{A}}

\newcommand{\F}{\mathbb{F}}

\newcommand{\spec}{\mathrm{Spec} \,}

\begin{document}

\begin{center}

{\large \sc A note on a generalization  \\ of the Hadamard quotient theorem}

\vspace{1cm}

Vesselin Dimitrov \\

\medskip

\texttt{vesselin.dimitrov@yale.edu}

\bigskip

\vspace{1cm}

\end{center}

\begin{abstract}
We consider a generalization of the ``Hadamard quotient theorem'' of Pourchet and van der Poorten. A particular case of our conjecture states that if $f := \sum_{n \geq 0} a(n)x^n$ and $g := \sum_{n \geq 0} b(n)x^n$ represent, respectively, an algebraic and a rational function over a global field $K$ such that $b(n) \neq 0$ for all $n$ and the coefficients of the power series $h := \sum_{n \geq 0} a(n)/b(n) x^n$ are contained in a finitely generated ring, then $h$ is algebraic. We prove this conjecture if either (i) $g$ has a pole of a strictly maximal absolute value at some place; or (ii) all poles of $g$ are simple, there is a positive density $\delta > 0$ of places which split completely in the field generated by the poles of $g$ and at which all $b(n)$ are invertible, and with $d := [K(t,f):K(f)]$ denoting the degree of the covering $f : C \to \proj^1$ defined by the algebraic function $f$, the local radii of convergence $R_v$ of $h$ at the places $v$ of $K$ satisfy $\sum_v \log^+{R_v^{-1}} \leq \delta \big/ 12d^4$.
\end{abstract}

\vspace{1cm}

\begin{center}
{\bf 1. Introduction}
\end{center}

{\bf 1.1. }  Let $K$ be any field. Consider the following two remarks. First, the coefficients of an algebraic power series in $K[[t]]$ are contained in a finitely generated ring. This is a classical theorem of Eisenstein, and can be seen for example by exhibiting an algebraic function as the diagonal of a rational function in two variables (or by basic principles of algebraic geometry). Second, if $f := \sum_{n \geq 0} a_nt^n \in K[[t]]$ is algebraic  and $g := \sum_{n \geq 0} b_nt^n \in K[[t]]$ is rational, then the {\it Hadamard product} $f * g := \sum_{n \geq 0} a_nb_nt^n \in K[[t]]$ is algebraic. More precisely, if $C$ is an algebraic curve over $K$,  $P \in C(K)$ a $K$-rational point, and $t \in K(C)$ a local parameter at $P$,
  the Hadamard product of the $P$-adic Taylor series in $K[[t]]$ of a rational function $f \in K(C)$ on $C$ with the Taylor expansion of a rational function in $K(t)$ is again the Taylor expansion around $P$ of a function in $K(C)$.

 The following converse to these observations is a generalization, apparently not considered before, of the so-called ``Hadamard quotient theorem'' proved by A. J. van der Poorten after a sketch~\cite{pourchet} by Y. Pourchet.  See~\cite{rumely}, and also its generalization by P. Corvaja and U. Zannier~\cite{corvaja}.

 \medskip

 {\bf Conjecture. } {\it Let $f := \sum_{n \geq 0} a_n t^n$ be algebraic  and $g := \sum_{n \geq 0} b_n t^n$ be rational, with $b_n \neq 0$ for all $n$. Assume the set $\{c_n := a_n/b_n \mid n \geq 0\}$ is contained in a finitely generated ring. Then $h := \sum_{n \geq 0} c_nt^n$ is algebraic. }

 \medskip

When $g = 1/(1-t)^2$ (hence $\mathrm{char}(K) = 0$), this is a generalization of a classical theorem of G. P\'olya, proved by Y. Andr\'e (see VIII 1.3 in~\cite{andreG}): an integral of an algebraic function is algebraic if its coefficients lie in a finitely generated ring. We remark that the latter result may also be seen as a very special case of a  conjecture of A. Ogus on the fullness of the de Rham realization of motives enriched with the crystalline Frobenii at the primes of a number field; see (2.3) in~\cite{ogus} for this relationship, and
 section 7.4 of Andr\'e's monograph on motives~\cite{andremotifs} for a discussion of Ogus' general conjecture.  However, Hadamard quotients of irrational algebraic functions by rational functions other than $1/(1-t)^2$ do not appear to have been considered in the literature.

\medskip

 {\bf 1.2. } When $K$ has characteristic $0$ and $f \in K(t)$ is also rational, Conjecture~1.1 is the content of the classical Hadamard quotient theorem as proved Pourchet and by van der Poorten, and generalized by Corvaja and Zannier to cover the weaker assumption that there is a finitely generated ring containing an infinite set of ratios $a_n/b_n$. (In doing this, one must be careful to exclude examples like $(2^n-2)/n$, much as in our Theorem 1.4 below; see~\cite{corvaja} for the precise statement.) The Corvaja-Zannier proof is based on Schmidt's Subspace theorem, and certainly does not carry over to positive characteristic. Van der Poorten's proof, which is based on a $p$-adic extrapolation procedure making essential use of the mixed characteristic of $\Q_p$, is also specific to number fields. The conjecture is therefore open even in the case $K = \F_q(x)$ and $f,g \in K(t)$.

\medskip

{\bf 1.3. } Conjecture~1.1 is obvious when $K$ is a finite field, and can be reduced by a standard specialization argument to the case of a global field: a number field (the mixed characteristic case), or the function field of a geometrically irreducible algebraic curve over a finite field (the positive equicharacteristic case). We therefore restrict our attention from now onwards to the case that $K$ is a global field.

Let $B$ be either the spectrum of the full ring of integers $O_K$ in a number field $K$; or a geometrically irreducible projective curve over a finite field $\mathbb{F}_q$. In the latter case, $K$ will denote the field of rational functions on $B$. For $S$ a finite set of places (closed points of $B$, plus all the archimedean places in the number field case), denote $O_{K,S}$ the ring of $S$-integers, or the rational functions on $B$ regular outside $S$. For $a \in K$ and $s \in B$ a closed point, let $\mathrm{ord}_s(a)$ denote the valuation of $a$ in $s$. For $\mathbf{x} = (x_0:x_1:\cdots:x_n) \in \proj^n(K)$ we consider the standard Weil height
\begin{eqnarray*}
h_K(\mathbf{x}) := -\sum_{s \in B} \min_j \,  \ord_s(x_j)  \cdot \log{|k(s)|}.
\end{eqnarray*}
Here $k(s)$ denotes the  residue field at $s$. In the geometric case,  the sum is only over the closed points $s$ of $B$. In the number field case, the term at the generic point of $\spec{O_K}$ is to be interpreted as the archimedean contribution, namely, the sum
$$
\sum_{\sigma : K \hookrightarrow \C} \max_j \, \log|\sigma(x_j)|,
$$
to be taken without the minus sign.

For $v \in B$ a closed point, consider the $v$-adic absolute value $|\cdot|_v$ defined by the usual normalization $\log{|\cdot|_v} := -\mathrm{ord}_v(\cdot) \cdot \log{|k(v)|}$. If $K$ is a number field and $v$ is a place corresponding to a complex embedding $\sigma : K \hookrightarrow \C$, we let $|\cdot|_v := |\sigma(\cdot)|^{n_v}$, where $n_v := [K_v:\R] \in \{1,2\}$. In those terms the height in rewrites in the more familiar form $h_K(\mathbf{x}) = \sum_{v} \max_j \log{|x_j|_v}$.

We also write $h_K(x):=h_K(1:x)$, viewing $x \in K$ as the point $(1:x)$ in $\proj^1(K)$.  For a polynomial $F \in K[\mathbf{x}]$ in several variables over $K$, we write $h_K(F)$ for the height of its set of coefficients, viewed as a point in a projective space. For $f = \sum_{n \geq 0} a_nt^n \in K[[t]]$, let $f_{/N} := \sum_{n=0}^N a_nt^n \in K[t]$ be the polynomial truncation modulo $t^{N+1}$. Following Bombieri in~\cite{bombieriG}, define the {\it height} of $f$ to be $h_K(f) := \limsup_{N} \frac{1}{N}h_K(f_{/N})$ as $N \to \infty$.

A set $S \subset B$ of closed points of $B$ is said to have {\it density $\delta$} if
$$
\lim_{n \to +\infty} \frac{\sum_{s \in S: \, |k(s)| \leq n} \log{|k(s)|}}{\sum_{s \in B : \, |k(s)| \leq n} \log{|k(s)|}} = \delta.
$$
 Finally, we fix an algebraic closure $\bar{K}/K$ and recall that the set of {\it poles} of an algebraic function over $K(t)$ consists of the zeros in $\bar{K}$ of the leading coefficient $p_0 \in K[t]$ of a minimal algebraic relation $\sum_{i=0}^d p_i f^{d-i} = 0$ with $p_i \in K[t]$.

\medskip

{\bf 1.4. } We will prove slightly more than Conjecture~1.1 in the case that the ``Hadamard denominator'' $\sum_{n \geq 0} b_nt^n \in K(t)$  has a dominant pole at some place.

\medskip

{\bf Theorem. } {\it Suppose that the Hadamard quotient $h := \sum_{n \geq 0} \frac{a_n}{b_n}t^n$ of an algebraic $f := \sum_{n \geq 0} a_n t^n \in K[[t]]$ by a rational $g := \sum_{n \geq 0} b_n t^n \in K[[t]]$ power series, where all $b_n \neq 0$, has finite height: $h_K(h) < +\infty$.

Assume there is a place $v_0$ of the splitting field of $g \in K(t)$ at which $g$ has a unique pole of a maximum  absolute value, possibly with a multiplicity. Then there is a non-zero $r \in K[1/(1-t)]$  such that the Hadamard product $h*r$ is algebraic. If moreover the dominant pole of $g$ at $v_0$ is a simple one, then $h$ is itself algebraic.
}

\medskip

{\bf 1.5. }
The multiplication by a power series in $K[1/(1-t)]$   really is necessary in view of the example of $-\log(1-t) = \sum_{n \geq 0}  t^n/n \in \Q[[t]]$, which has finite height $e < + \infty$ and is the Hadamard quotient of $1/(1-t)$ by $1/(1-t)^2$, but is not algebraic. When $r = \sum_{n \geq 0} p(n)t^n \in \Q(t)$ with $p \in \Q[x]$ a polynomial which splits in linear factors over $\Q$, it follows easily from Andr\'e's result mentioned in 1.1 that the algebraicity of $h$ itself is implied by Theorem 1.4 together with the original integrality assumption of Conjecture 1.1, stronger than the finiteness of the height, namely the finite generation of the ring $\Z[ a_n/b_n \mid n \geq 0]$. In this case, assuming $p(n) \neq 0$ for all $n$, $\sum_{n \geq 0} t^n/p(n) \in \Q[[t]]$ indeed has finite height, being a $\Q(t)$-linear combination of polylogarithmic series; and we conjecture that $r = \sum_{n \geq 0} p(n)t^n$ in Theorem~1.4 may always be taken to be of this form, with $p \in \Q[x]$ completely split over $\Q$. This is subsumed by Conjecture 1.8 (iv) below.

\medskip

{\bf 1.6. }
The following condition, intermediate between the finiteness of the height of $h(t) = \sum_{n \geq 0} c_nt^n \in K[[t]]$ considered in Theorem 1.4 and the finite generation of the ring $\Z[c_n \mid n \geq 0]$ considered in Conjecture 1.1,
arises naturally in this type of questions, and allows for a clean conjectural generalization of 1.1 in the number field case which does not involve the corrective factor $r \in K[1/(1-t)]$ in Theorem~1.4. It is suggested to us by the paper~\cite{bostacl} by J.-B. Bost and A. Chambert-Loir, from which we borrow the terminology of {\it $A$-analyticity} (itself borrowed from a suggestion of Y. Manin from a 1984 Arbeitstagung).

\medskip

{\bf Definition. } {\it A power series $f \in K[[t]]$ is \emph{$A$-analytic} if there exists a collection of numbers $r_v \in (0,1]$, indexed by the closed points (finite places) $v \in B$, such that $\prod_{v \in B} r_v > 0$ and for each $v$, the $v$-adic radius of convergence of $f$ is at least $r_v$, and $f$ maps the disk $|t|_v < r_v$ to itself.}

\medskip

For example, if there is a finite set of places $S$ such that $f \in O_{K,S}[[t]]$ has $S$-integral coefficients and positive radii of convergence at $v \in S$, we may take $r_v = 1$ for all $v \notin S$, hence such power series are $A$-analytic. An $A$-analytic power series has finite height, but the reverse does not hold.
In the number field case, for instance, $\log(1+t) =\sum_{n \geq 1} (-1)^{n+1}t^n/n$ is not $A$-analytic: while the radius of convergence is $1$ at every place,  $r_p$ may not exceed $p^{-1/(p-1)}$ for every prime $p$, and the sum $\sum_p \frac{\log{p}}{p-1}$ over the primes diverges.

This notion, which is reminiscent of the Bombieri conditions from the theory of $G$-functions, arose originally in Bost's work~\cite{bostfoliations} on a non-linear generalization by Ekedahl, Shepherd-Barron, and Taylor of the Grothendieck-Katz $p$-curvature conjecture; see condition (3.1) of Theorem~3.4 in {\it loc. cit.}, where the general notion of $v$-adic size of a formal subscheme is defined in 3.1.1  {\it loc. cit.} We note that when $h \in K[[t]]$ is $A$-analytic, our ``dominant pole case'' Theorem~1.4.   is in fact an easy consequence of the theorem of Bost just cited (in the case that the place $v_0$ is archimedean), and its generalization Theorem~6.1 in~\cite{bosteval} (in the general case over a number field); Bost's arguments carry through verbatim to also cover the function field case. The key remark  is that, by a geometric series expansion, the Hadamard inverse of a rational function $g \in K(t)$ having a dominant pole at the place $v$ is sufficiently well-approximated $v$-adically by a rational function on $\proj^1$ to imply Bost's ``condition $L_v$'' (5.2 in~\cite{bosteval}) for the formal graph in $C \times \proj^1$ of the Hadamard quotient of a rational function on $C$ by such a $g \in K(t)$. We provide a direct proof  which remains valid under the weaker hypothesis $h_K(h(t)) < +\infty$. It combines a much older idea of D. Cantor~\cite{cantor}, which it extends to higher genus, with an analytic technique of Bost from~\cite{bostfoliations,bostgerms} based on the Poincar\'e-Lelong equation.

\medskip

We complement Theorem 1.4 with the following rather different  case of Conjecture 1.1, which is valid under the $A$-analyticity condition and which, though rather restrictive, does include cases not covered by the dominant pole condition.

\medskip

{\bf Theorem 1.7. } {\it Let $f = \sum_{n \geq 0} a_n t^n \in K[[t]]$ be an algebraic power series, and consider $d := [K(t,f):K(f)]$ the degree of the induced branched cover $f : C \to \proj_K^1$ of the regular and projective model of the field extension $K(t,f)/K(f)$. Let $g = \sum_{n \geq 0} b_nt^n \in K[[t]] \cap K(t)$ be a rational power series with only simple poles.

(i) Assume there is a positive density $\delta > 0$ of closed points of $B$ (finite places of $K$) which split completely in the extension of $K$ generated by the poles of $g$, and at which all coefficients $b_n$ are units.  If the Hadamard quotient $h := \sum_{n \geq 0} \frac{a_n}{b_n} t^n \in K[[t]]$ is $A$-analytic and its $v$-adic radii of convergence $R_v > 0$ satisfy the condition
$$
\sum_v \log^+{R_v^{-1}} \leq \delta/12d^4,
$$
then $h$ is algebraic.

(ii) If moreover $f \in K(t)$ is also rational with only simple poles,  and if there is a positive density $\delta > 0$ of closed points of $B$ which split completely in the extension of $K$ generated by the poles of $g$ \emph{and} $f$, and at which all $b_n$ are units, then $h$ is rational as soon as its height satisfies
$$
h_K(h) < \delta/12.
$$}

\medskip

For example, though this case is covered by the dominant pole condition of Theorem 1.4, we note a result of  H. Hasse~\cite{hasse} stating that if $a \geq 3$ is a square-free integer (resp. $a = 2$), the rational primes not dividing any term of the sequence $a^n + 1$ of coefficients of $(2-(a+1)x)\big/(1-x)(1-ax) \in \Z[[x]] \cap \Q(x)$ have density equal to $1/3$ (resp. $7/24$).

In view of the Corvaja-Zannier theorem, condition (ii) is of interest only in the positive equicharacteristic (``function field'') case. Let us also note that, outside of degenerate situations, there probably do not exist any non-trivial algebraic (resp., rational) functions $h$ satisfying the condition required for the conclusion in (i) (resp., (ii)). This does not mean that Theorem 1.7 is vacuous, but rather that it should be read in the contrapositive sense, as a lower bound outside of degenerate situations on the height of an Hadamard quotient.

\medskip

{\bf 1.8. } We close the introduction by suggesting that much more than Conjecture 1.1 could be true.   Extending the standard terminology to include the positive characteristic case, we say that $f \in K[[t]]$ is a {\it $G$-series} if $f$ is $D$-finite (meaning it satisfies a linear differential equation with polynomial coefficients) and $h_K(f) < \infty$. The standard source of examples are the power series expansions of a variation of periods in a family of algebraic varieties over $K$. Those satisfy a linear differential equation (Picard-Fuchs) coming from the finite dimensionality of cohomology. Let us say here that the $G$-series which ``come from geometry'' in this way are of motivic origin; see Y. Andr\'e's book~\cite{andreG} for more details and a precise definition of $G$-functions ``coming from geometry.'' A fascinating conjecture of E. Bombieri and B. Dwork suggests that all $G$-series are motivic in this sense.

A particular example of motivic $G$-series are the diagonals of rational functions in several variables, see~\cite{christol,andreG}. The latter have been conjectured in~\cite{christol} by G. Christol to be exactly the $G$-series whose coefficients are contained in a finitely generated ring. We suggest the following generalization of the Hadamard quotient theorem, a motivation for which is that the classes of $G$-series, motivic $G$-series, and diagonals of rational functions are each closed under Hadamard product $\big(\sum_{n \geq 0} a_n t^n \big) * \big(\sum_{n \geq 0} b_n t^n \big) \mapsto \big(\sum_{n \geq 0} a_nb_n t^n \big)$.
This is clear for $G$-series, and for the proof in the other two subclasses we again refer to~\cite{andreG}. The complete {\it diagonal} of a power series $\sum_{\mathbf{n}} a_{\mathbf{n}}\mathbf{t}^{\mathbf{n}} \in K[[t_1,\ldots,t_r]]$ in several variables is the single variable power series $\sum_{n} a_{n,\ldots,n}t^n \in K[[t]]$. Regarding part (vii) below we recall that algebraic functions in positive characteristic are closed under the Hadamard product; two different proofs of this fact are contained in~\cite{furstenberg} and~\cite{deligne}.

\medskip

{\bf Conjecture. } {\it \begin{itemize}
 \item[(i)] If the Hadamard quotient of two $G$-series has finite height, then again it is a $G$-series. That is, the Hadamard quotient is $D$-finite if it has finite height;
     \item[(ii)] If the Hadamard quotient of two motivic $G$-series has finite height, then again it is a motivic $G$-series. That is, it satisfies a differential equation of Picard-Fuchs type;
\item[(iii)] If the Hadamard quotient of an algebraic by a rational power series is $A$-analytic, then it is algebraic;
\item[(iv)] If the Hadamard quotient of an algebraic by a rational power series has finite height, then there exists a non-zero polynomial $p \in \Q[x]$ which splits into linear factors over $\Q$, such that, denoting $r(t) := \sum_{n \geq 0} p(n)t^n \in K[1/(1-t)]$, the Hadamard product $h * r$ is algebraic.
         \item[(v)] If the Hadamard quotient of diagonals of rational functions in several variables is $A$-analytic, then again it is the diagonal of a rational function in several variables.
\item[(vi)] If the Hadamard quotient $h$ of diagonals of rational functions in several variables has finite height, then there exists an $r(t)$ as in (iv) such that $h * r$ is again the diagonal of a rational function in several variables.
             \item[(vii)] Assume $\mathrm{char}(K) > 0$. If the Hadamard quotient of two algebraic power series has finite height, then again it is algebraic.
\end{itemize}}

\medskip

Though we will not consider these questions in the present paper, we remark that the method of proof of Theorem 1.4 extends to confirm Conjecture 1.8 in many cases where there exists a place $v_0$ at which the coefficients $b_n$ of the ``Hadamard denominator'' satisfy the asymptotic growth condition $\log{|b_n|_{v_0}} \sim cn$. This is our motivation for the general statement we propose here.

Finally, we note that one may consider related questions in the spirit of Pisot's $d$-th root conjecture, proved by Zannier in~\cite{zannier} for the case of Hadamard $d$-th root from rational functions over a number field. For instance, if $f = \sum_{n \geq 0} a_nt^n \in K[[t]]$  and $P \in K[x,y]$ are such that $f$ is a $G$-function and $P(a_n,y) = 0$ has for each $n$ a $K$-rational solution $y \in K$, does it follow that there exists a $G$-function $\sum_{n \geq 0} b_nt^n \in K[[t]]$ with $P(a_n,b_n) = 0$ for each $n$?

\bigskip

\begin{center}
{\bf 2. Proof of Theorem 1.4}
\end{center}

\bigskip

Upon replacing $K$ by a finite extension, we may assume that $g \in K(t)$ splits over $K$. By assumption, recalling the explicit shape of coefficients  of rational power series as confluent power sums, there is a place $v_0$ of $K$ at which the coefficients $g = \sum_{n \geq 0} b_nt^n$ are of the form
\begin{equation} \label{dominant}
b_n = \beta^n \Big( p(n) + \sum_{i=1}^k q_i(n) \gamma_i^n  \Big), \quad p \in K[x] \setminus \{0\}, \, q_i \in K[x],  \textrm{ and all } |\gamma_i|_{v_0} < 1.
\end{equation}
If moreover the dominant pole at $v_0$ is a simple one, we may take $p = 1$. In the positive characteristic case it is important to stress that $p$ does not vanish identically on $\Z$. 

Since the power series $f \in K[[t]]$ is algebraic, it has a strictly positive radius of convergence $\rho_0 > 0$ at the place $v_0$. On the other hand, by the geometric series expansion and the strict inequalities $|\gamma_i|_{v_0} < 1$ in (\ref{dominant}),  we may write
\begin{equation} \label{approx}
p(n)^l\beta^n/b_n = u_n + r_n,
\end{equation}
where $l \in \N$ and the power series $u := \sum_{n \geq 0} u_nt^n \in K(t)$ is rational, while the radius of convergence $R_0'$ of the remainder $r := \sum_{n \geq 0} r_nt^n \in K[[t]]$ at the place $v_0$ can be as large as desired.
 Then the Hadamard product $z := f * r \in K[[t]]$ converges on the disc of radius $R_0 := \rho_0R_0'$, which in turn may be taken to be arbitrarily large.

Let $C/K$ be the smooth projective model of the finite field extension $K(t,f)/K(f)$, and $C \to \proj_K^1$ the corresponding (branched) covering. Thus $t \in K(C)$ is a non-constant rational function on $C$. By design, since $f \in K[[t]]$ implies $K(C) \subset K((t))$, the  function $t \in K(C)$ has only simple zeros.
 Therefore the divisor of zeros of $t$ is a finite reduced non-empty  subscheme $Z$ of $C/K$, and $t$ is a local parameter at all points in $|Z|$. We  can view the Hadamard quotient $h \in K[[t]]$ as a formal germ in $\widehat{\Oz_{C,Z}}$, the formal completion  of $C$ at  $Z$.

Note that the Hadamard product of $f = \sum_{n \geq 0} a_nt^n$ with a power series in $K(t) \cap K[[t]]$ is again in $K(C)$; this follows from the explicit description of the coefficients of rational power series as confluent power sums, which reduces the verification to the case $\sum_{n \geq 0} a_n \binom{n+i}{i} \alpha^n t^n = \alpha^{-i} \partial^i f(\alpha t)/\partial t^i$. Therefore $f * u \in K[[t]]$ is the expansion around $Z$ of a rational function $y \in K(C)$ on $C$, and~(\ref{approx}) expresses (a derivative of) the Hadamard quotient $h$ as
\begin{equation}  \label{hp}
p^l(t \,\partial/\partial t)h(\beta t) = \sum_{n \geq 0} p^l(n)\frac{\beta^n a_n}{b_n} t^n = y + z, \quad y \in K(C), \, z \textrm{ convergent on } |t|_{v_0} < R_0.
\end{equation}

Let $D$ be the polar divisor of $t$, and
consider the  ample line bundle $L := \Oz(D)$ on $C$ and its canonical global section $s_0 := 1_D$ with divisor $D$. By the Riemann-Roch theorem, there exists a $k < +\infty$ and a non-zero global section $s \in \Gamma(L^{\otimes k})$ such that the section $y \cdot s$ of $L^{\otimes k}$ is regular outside $D$. Letting
$$
U := \{  |t|_{v_0} < R_0 \} \subset C_{v_0}^{\mathrm{an}} \setminus |D|,
$$
we therefore have $y \cdot s, \, z \cdot s \in \Gamma(L^{\otimes k}, U)$. {\it We claim that the conditions $h_K(h) < +\infty$ and $R_0 \gg 0$ imply that the formal germ $F(t) := \big(p^l(t \, \partial/\partial t) h(\beta t)\big) \cdot s/s_0^k = (y+z) \cdot s/s_0^k \in \widehat{\Oz_{C,Z}}$ is a polynomial in $t$, hence in particular a rational function on $C$.} As $\beta \neq 0$ and $s/s_0^k \in K(C)$ is also a non-zero rational function on $C$, this will prove the algebraicity of $p^l(t \, \partial/ \partial t)h$, which is the Hadamard product of $h$ with a power series in $K[1/(1-t)]$ having non-zero coefficients. If as in the last clause of Theorem 1.4 the dominant pole at $v_0$ is a simple one, we have $p = 1$, and it will follow in this case that the Hadamard quotient $h$ is itself algebraic. Assuming to the contrary that the Taylor series of $F(t) \in K[[t]]$ is infinite, we will derive a contradiction by the product formula from the existence of arbitrarily large  values of $N < +\infty$ for which $t^N$ appears in $F$ with a non-zero coefficient $e_N \in K \setminus \{0\}$.

We thus consider a large $N$ such that $e_Nt^N$ appears in the series $F(t) \in K[[t]]$ with a non-zero coefficient $e_N \in K \setminus \{0\}$. The product formula then gives
\begin{equation} \label{prod}
\sum_v \log{|e_N|_v} = 0,
\end{equation}
where the sum is over all places $v$ of $K$.  Since the sum of the contributions to~(\ref{prod}) over any subset of the places is bounded   by the height $h_K(e_N)$, we have
\begin{equation} \label{este}
-\log{|e_N|_{v_0}} = \sum_{v \neq v_0} \log{|e_N|_v} \leq h_K(e_N).
\end{equation}
To estimate the height, we use the bound (see Prop. 1.5.15 in~\cite{bible})
\begin{equation} \label{easybound}
h_K(\alpha_1 + \cdots + \alpha_r) \leq \varepsilon \log{r} + \sum_v \max_j \log{|\alpha_j|_v},
\end{equation}
applied to the sum defining $e_N$ as the coefficient of $t^N$ in $(s/s_0^k)_{/N} \cdot p^l(t \, \partial/\partial t)h(\beta t)_{/N} = (s/s_0^k)_{/N} \cdot \sum_{n = 0}^N p(n)^lc_n\beta^nt^n$. (Here, $\varepsilon$ denotes $[K:\Q]$, if $K$ is a number field, and $0$ otherwise.) We obtain
\begin{equation} \label{esteh}
h_K(e_N) \leq (\varepsilon + c) \log{N} + h_K(h_{/N}) + h_K(\beta) + h_K\big( (s/s_0^k)_{/N} \big),
\end{equation}
where for $c \gg 0$ the term $c \log{N}$ takes care of the polynomial factor $p(n)$.
Our assumption that the Hadamard quotient $h$ has finite height means that $h_K(h_{/N}) = O(N)$, while the Taylor series of $s/s_0^k$, being  algebraic, certainly has finite height by Eisenstein's theorem. From (\ref{este}) and (\ref{esteh}) we derive the upper bound
\begin{equation} \label{upper}
-\log{|e_N|_{v_0}} \leq A_0 \cdot N,
\end{equation}
where  $A_0 < +\infty$ is a constant independent of $N$.

On the other hand, as noted above, the parameter $R_0$ can be taken to be as large as desired. Take $R_0 > e^{A_0}$, and fix it. (In the ultrametric case it will be convenient to assume additionally that $R_0 \in |K_{v_0}|_{v_0}$.) We will next show that $\log{|e_N|_{v_0}} \leq -N \log{R_0} + O(1)$ as $N \to \infty$, which together with (\ref{upper}) will complete the proof by contradiction.

The truncation $w := F_{/N-1}$ modulo $t^N$ is  the unique polynomial of degree $< N$ in $K[t]$ for which
\begin{equation} \label{coeff}
(y+z) \cdot s/s_0^k - w = e_Nt^N + o(t^N).
\end{equation}
In particular, the section
\begin{equation} \label{given}
s' := y \cdot s + z \cdot s - w \cdot s_0^k \in \Gamma(L^{\otimes k}, U)
\end{equation}
vanishes to order $N$ along $Z$.

We divide the majorization of $\log{|e_N|_{v_0}}$ into two alternatives, whether the place $v_0$ is Archimedean or ultrametric. In the former case we borrow a technique from Bost's paper~\cite{bostgerms} based on the Poincar\'e-Lelong equation (see the proof of Prop 3.6 in {\it loc. cit.}); this idea germinates in fact from Bost's earlier work~\cite{bostfoliations} in relation to analytic parametrizations by complex manifolds large enough to have Liouville's property that a bounded harmonic function is constant. Its ultrametric counterpart appears in the proof of Prop. 5.14 of~\cite{bostacl}.

\medskip

{\it Archimedean case.} Suppose $v_0$ is Archimedean, corresponding to a complex embedding $\iota : K \hookrightarrow \C$.
 Choose any $C^{\infty}$ hermitian metric $\|\cdot\|$ in the line bundle $L_{v_0} := (L \otimes_K \C)^{\mathrm{an}}$ on the closed Riemann surface $C_{v_0}^{\mathrm{an}}$; it induces a metric in all tensor powers $L_{v_0}^{\otimes n}$, which we will continue to denote by $\|\cdot\|$. As in the proof of Prop. 3.6 in Bost's paper~\cite{bostgerms}, consider any $C^{\infty}$ function $\psi : C_{v_0}^{\mathrm{an}} \setminus |D| \to \R$ such that $\psi|_Z = 0$ and
\begin{equation} \label{subpotential}
\frac{\sqrt{-1}}{\pi} \partial \bar{\partial} \psi \geq k c_1(L) = c_1(L^{\otimes k}),
\end{equation}
as a point-wise inequality of $(1,1)$-forms on the Riemann surface $S := C_{v_0}^{\mathrm{an}} \setminus |D|$, where $c_1(L) := \frac{\sqrt{-1}}{\pi} \partial\bar{\partial} \log{\|s_0\|} \in A^{1,1}(S)$ is the Chern form of $(L,\|\cdot\|)$. For example,  we note following Bost that for any two holomorphic functions $h_1,h_2$ on $C \setminus |D|$ which vanish along $Z$ and have disjoint ramification divisors,  the choice $\psi := A(|h_1|_{v_0}^2 + |h_2|_{v_0}^2)$ is admissible for a large enough constant $A \gg 0$. As $Z$ is reduced, such $h_1,h_2$ clearly exist.

Since $s'$ is a section of $L^{\otimes k}$, we have the Poincar\'e-Lelong distributional equation
\begin{equation} \label{residue}
\frac{\sqrt{-1}}{\pi} \partial \bar{\partial} \log{\|s'\|} = \delta_{\mathrm{div}(s')} - c_1(L^{\otimes k})
\end{equation}
in the sense of currents on $S$. Since by construction the divisor of  $s'$ on $U \subset S$ dominates $N \cdot Z$, we can retain from the sum of~(\ref{subpotential}) and~(\ref{residue}) the inequality
\begin{equation} \label{current}
\frac{\sqrt{-1}}{\pi} \partial \bar{\partial} \big(  \log{\|s'\|} + \psi  \big) \geq N \delta_Z =  N\frac{\sqrt{-1}}{\pi} \partial \bar{\partial} \log{|t|_{v_0}}
\end{equation}
of  $(1,1)$-currents on $U$. Here, for the second equality, we have used that $Z$ is the divisor on $U$ of the function $t$, and  again the Poincar\'e-Lelong equation.

 This means that the function
$$
 \log{\|s'\|} + \psi + N\log{|R_0/t|_{v_0}} = \log{\frac{\|s'\|}{|t|_{v_0}^N}} + \psi + N\log{R_0}
$$
from $S$ to $\{-\infty\} \cup \R$ is subharmonic on $U = \{ |t|_{v_0} < R_0 \}$. Since $\log{|R_0/t|_{v_0}} = 0$ on the boundary $\partial U$, the value of this subharmonic function is bounded on $\partial U$ by $B_1 := \sup_{\partial U} \big( \log{\|s'\|} + \psi \big) < +\infty$, and the maximum principle implies that this last quantity bounds the value at a point $P \in |Z|$ (which we choose and fix). Since $\psi|_Z = 0$ by construction, the conclusion reads:
\begin{equation} \label{cauchy}
\lim_{p \to P} \log \frac{\|s'(p)\|}{|t(p)|^N} \leq -N\log{R_0} + B_1.
\end{equation}
By~(\ref{coeff}) and (\ref{given}), there exists a constant $B_2 < +\infty$, depending on $k$ and the choice of metric $\|\cdot\|$ but not on $N$, such that the the left-hand side of~(\ref{cauchy}) is at least $\log{|e_N|_{v_0}} - B_2$. Thus, with $B_3 := B_1 + B_2$, we have obtained the required majorization
\begin{equation} \label{dominantestimate}
\log{|e_N|_{v_0}} \leq -N\log{R_0} + B_3.
\end{equation}

\medskip

{\it Ultrametric case. } In the case that $v_0$ is ultrametric the argument is easier, and owes to the proof of Prop. 5.14 (comparison of canonical and capacitary metrics) in Bost and Chambert-Loir's paper~\cite{bostacl}. Without loss of generality we may assume that $\log{R_0}$ is in the value group of the discrete valuation $|\cdot|_{v_0}$; write $R_0 = |\xi|_{v_0}$ with $\xi \in K_{v_0}^{\times}$.
Consider the rigid $K_{v_0}$-analytic curve
$$
U := \{ |t|_{v_0} < R_0 \} \subset C_{v_0}^{\mathrm{an}}
$$
in the $v_0$-adic analytification of $C$; it is the complement in $C_{v_0}^{\mathrm{an}}$ of the affinoid domain $U^{\mathrm{c}} := C_{v_0}^{\mathrm{an}} \setminus U = \{ |\xi/t|_{v_0} \leq 1 \}$, and its {\it boundary} is the affinoid $\partial U := \{ |\xi/t|_{v_0} = 1 \}$. The affinoid function $s'/s_0^k \cdot (\xi/t)^N$ on $U$ is analytic by construction, and its value at each point in $|Z|$ equals $e_N\xi^N$. By the ultrametric maximum principle (see Prop. B.1 in~\cite{bostacl}) it follows that
$$
\log{|e_N\xi^N|_{v_0}} \leq \sup_{\partial U} |s'/s_0^k \cdot (\xi/t)^N|_{v_0} =  \sup_{\partial U} |s'/s_0^k|,
$$
from which we retain, with $B_3 := \sup_{\partial U} |s'/s_0^k|$, the required majorization
\begin{equation} \label{dominantestbis}
\log{|e_N|_{v_0}} \leq -N \log{|\xi|_{v_0}} + B_3 = -N\log{R_0} + B_3.
\end{equation}

\medskip

In either case, as noted above, the minorization obtained contradicts~(\ref{upper}) as soon as $\log{R_0} > A_0$ and $N \gg_{R_0} 0$. Thus the bounds~(\ref{dominantestimate},\ref{dominantestbis}) yield the contradiction completing the proof of the theorem. \proofend

\bigskip

\begin{center}
{\bf 3. An algebraicity criterion}
\end{center}

\medskip

We shall deduce Theorem~1.7 from an algebraicity criterion on a formal power series in $K[[t]]$.  In an appendix to this paper we propose an optimal conjectural generalization of this criterion which would generalize at the same time an old conjecture of I. Ruzsa and the classical rationality criterion of Borel-Dwork.

As in 1.3, $K$ will denote a global field. We reserve the following notation for a uniform treatment of the function field and number field cases. Let $\varepsilon$ be $0$ in the former case and $[K:\Q]$ in the latter case. In the former case let $G :=g\log{q}$ where $g$ is the genus of the regular projective curve $B_{/\mathbb{F}_q}$, while in the later case  write $G := \frac{1}{2} \log{|D_{K/\Q}|}$, where $D_{K/\Q}$ is the discriminant of the number field $K$.

\medskip

{\bf 3.1. Siegel's lemma.}
 The relevant Siegel lemmas are as follows.
 \begin{quote}
  {\it  If $A$ is an $M \times N$ matrix, $N > M$, with coefficients  $a_{ij} \in K$ satisfying $h_K(a_{ij}) \leq H$,  then $A.\mathbf{x} = 0$ admits a non-zero solution $\mathbf{x}$ of height satisfying
 $$
 h_K(\mathbf{x})   \leq \frac{\varepsilon}{2} \log{N} + G  + \frac{MH}{N-M}.
 $$
  }
 \end{quote}
 In the number field case this follows from the theorem of Bombieri and Vaaler, for which we refer to section 2.9 in~\cite{bible}. In the function field case this is a particular instance of R. Thunder's Siegel lemma from~\cite{thunder}.

 \medskip

Recall from 1.3 the notation $f_{/n}$ for the polynomial truncations of a power series. For polynomials $P_1,\ldots,P_m \in K[t]$, we use the notation $h_K(P_1,\ldots,P_m)$ for the height of the set of all coefficients of all $P_i$, viewed as a point in a projective space.

\medskip

 {\bf 3.2. An algebraicity criterion. }
 {\it Let $f \in K[[t]]$.
  Fix a positive integer $r$. For each closed point $s \in B$ at which $f$ can be reduced,  let $h_s = h_s^{(r)}(f) \in \N_0 \cup \{+\infty\}$ be the maximum of the degrees of the coefficients of a minimal algebraic relation of $f_s$ of degree $\leq r$ over $K(t)$, with the convention that $h_s = +\infty$ if no such relation exists.
  If
\begin{eqnarray} \label{criterion}
    \liminf_n \Big\{ \frac{1}{n} \sum_{s: \, h_s < n/4r^2} \log{|k(s)|} \Big\}  > (2r+1)  \limsup_n \Big\{ \frac{1}{n} h_K(1,f_{/n},\ldots, (f^r)_{/n}) \Big\},
\end{eqnarray}
 then $f$ is algebraic of degree at most $r$ over $K(t)$.}

 \medskip

{\bf Remark. } The constants $1/4r^2$ and $2r+1$ can be improved in the proof below; we do not do this here in order to keep the qualitative statement as simple as possible.

\medskip

 {\it Proof. }   Consider $L$ a large integer parameter. Applying Siegel's lemma~2.1 with
$$
M := 2r(r+1)L \textrm{ and }  N := (2r+1)(r+1)L.
$$
we find a polynomial $\Phi \in \Gamma \big(\proj_K^1 \times \proj_K^1, \Oz ( (2r+1)L,r)\big)$ of bi-degree $((2r+1)L,r)$  such that $\Phi(t,f(t))  \equiv 0 \mod{t^{2r(r+1)L}}$ and
\begin{eqnarray} \label{est2}
h_K(\Phi) \leq \frac{\varepsilon}{2}\log{N}  + 2r \cdot  h_K\big(1,f_{/M}, \ldots, (f^r)_{/M} \big) + G.
\end{eqnarray}
Assuming $\Phi(t,f(t))  \neq 0$, let $n \geq 2r(r+1)L$ be the order of vanishing at $t = 0$, and let $c \in K \setminus \{0\}$ be the coefficient of the first non-zero term $t^{n}$ in the expansion. By the product formula,  we have
\begin{equation} \label{prod2}
\sum_{\textrm{all } v} \log{|c|_v} = 0.
\end{equation}

The equation $\Phi  = 0$ cuts on $\proj^1 \times \proj^1$ a curve bi-degree $((2r+1)L,r)$ with respect to the canonical very ample line bundle $\Oz(1,1)$. At every closed point $s \in B$ the graph $\Gamma_s$ of the reduction $f_s$ is a curve of bi-degrees at most $(h_s,r)$ on $\proj_{k(s)}^1 \times \proj_{k(s)}^1$. In the fibre over $s \in B$, the intersection number of this curve with the curve $\{\Phi = 0\}$ does not exceed $r(2r+1)L+rh_s$, while the intersection index at the point $P := (0,f(0))$  equals $n$. Since the curve $\Gamma_s$ is smooth, in particular irreducible, it follows that $\Gamma_s$ is a component of the curve $\Phi \equiv 0 \mod{s}$ as soon as $r(2r+1)L+rh_s < n$. Since by construction $n \geq 2r(r+1)L$, the last inequality is implied by $h_s < n/2r(r+1) < n/4r^2$.

 At every closed point $s \in B$ with $h_s < n/4r^2$ we have by the previous paragraph that $\mathrm{ord}_s(c) > 0$, hence the contribution at $s$ to~(\ref{prod2}) is at least $\log{|k(s)|}$. The absolute value of the sum of the remaining contributions does not exceed $h_K(c)$, hence the product formula~(\ref{prod2}) yields the lower bound
\begin{equation}  \label{eqn}
h_K(c) \geq \sum_{s: \, h_s < n/4r^2} \log{|k(s)|}.
\end{equation}
For an upper bound we use~(\ref{est2}), the general bound~(\ref{easybound}) used also in the proof of Theorem 1.4, and the definition of $ct^n$ as a term in $\Phi(t,f(t)) \in K[[t]]$. We obtain:
\begin{eqnarray*}
h_K(c) \leq \varepsilon \log{N} + h_K(\Phi) +h_K(f_{/n},\ldots,(f^r)_{/n}) \\
\leq \frac{3\varepsilon}{2} \log{N}  + G +  (2r+1) h_K(1,f_{/n},\ldots,f_{/n}^r).
\end{eqnarray*}
As $L$ and hence $n$ can be taken arbitrarily large, the last bound combined with~(\ref{eqn}) contradicts the assumption of the theorem. The contradiction forces upon $f$  the degree-$r$ algebraic relation $\Phi(x,f(x))  = 0$. \proofend

\medskip

When $f \in K[[t]]$ is $A$-analytic,  the criterion 3.2 is complemented by the following evaluation of the right-hand side of  condition~(\ref{criterion}).

\medskip

{\bf Lemma 3.3. } {\it If $f \in K[[t]]$ is $A$-analytic and has for each place $v$ of $K$ a $v$-adic radius of convergence $R_v > 0$, then for each $r \in \N$,
$$
\limsup_n \Big\{ \frac{1}{n} h_K(1,f_{/n},\ldots, (f^r)_{/n}) \Big\} = h_K(f) = \sum_v \log^+{R_v^{-1}}.
$$}

\medskip

{\it Proof. } Let $r_v \in (0,1]$ be as in Definition 1.6, and
$$
C := \sum_v \log^+{r_v^{-1}} < +\infty.
$$
Each power $f^j$ has $v$-adic radius of convergence $R_v$ and maps the disk $|t|_v < r_v$ to itself. Therefore, the expansion
$$
f^j = \sum_{n \geq 0} a_n^{(j)}t^n \in K[[t]]
$$
has coefficients satisfying
\begin{equation} \label{dom}
\log{|a_n^{(j)}|_v}  \leq (n-1)\log{r_v^{-1}}, \quad \limsup_{n} \frac{1}{n}\log{|a_n^{(j)}|_v} = \log{R_v^{-1}}.
\end{equation}
Here, the first estimate holds for all finite places $v$, and the second equality is valid for all $v$. The latter implies in general that the required limit supremum is at least $\sum_v \log^+{R_v^{-1}}$, and we have to show that under the $A$-analyticity condition,
\begin{equation} \label{desired}
h_K(1,f_{/n},\ldots, (f^r)_{/n}) < n\sum_v \log^+{R_v^{-1}} + n\epsilon
\end{equation}
for $n \gg_{\epsilon} 0$ and any $\epsilon > 0$.

To this end, take a finite set $S$ of places of $K$, including all Archimedean ones, such that
$$
\sum_{v \notin S} \log^+{r_v^{-1}} < \epsilon/2,
$$
hence, by the first inequality in (\ref{dom}),
\begin{equation} \label{domi}
\sum_{v \notin S} \max_j \log^+{|a_n^{(j)}|_v} < n \epsilon/2.
\end{equation}
Since $S$ is finite, the second equality in (\ref{dom}) implies
\begin{equation} \label{domid}
\max_{0 \leq j \leq r} \log^+{|a_n^{(j)}|}_v < n\log^+{R_{v}^{-1}} + n\frac{\epsilon}{2|S|} \quad \textrm{ for all } v \in S,  \, n \gg_{S,r} 0.
\end{equation}
The desired bound~(\ref{desired}) then follows by summing  (\ref{domi}) with the estimates (\ref{domid}) for all $v \in S$. \proofend

\bigskip

\begin{center}
{\bf 4. Proof of Theorem 1.7}
\end{center}

\medskip

We first consider part (i). In the notation of 3.2, the assumption that $f$ defines a degree-$d$ branched covering $C \to \proj^1$ implies $h_s^{(d)}(f) \leq d$ for all closed points $s \in B$. Starting from this, the key  the observation,  which was also used in the proof of Theorem 1.4, is that the Hadamard product of an algebraic and a rational power series is algebraic. More precisely, by the explicit description of coefficients of rational functions as confluent power sums shows that over any field $k$, the Hadamard product of an element of $k[[t]]$ which satisfies an algebraic relation of degree $r$ over $k(t)$ with coefficients of degree at most $a$, and an element of $k[[t]] \cap k(t)$ of degree at most $b$, satisfies an algebraic relation of degree $r$ over $k(t)$ with coefficients of degree at most $ab$.

By assumption, there is a set $S$ of closed points of $B$ which has density $\delta$, and such that all coefficients $b(n)$ of the ``Hadamard denominator'' $g = \sum_{n \geq 0} b(n)t^n \in K[[t]] \cap K[t]$  have at each $v \in S$ a reduction $\tilde{b}_v(n)$ which at once lies in $k(v)^{\times}$  and has the form $\sum_{i=1}^k c_ib_i^n$ with $c_i,b_i \in k(v)$. Such a sequence $\tilde{b}_v(n)$ is periodic with pure period $|k(v)^{\times}| = |k(v)| - 1$, and so also is its inverse $1/\tilde{b}_v(n) \in k(v)^{\times}$. This shows that the Hadamard inverse of $g\mod{v}$ has degree less than $|k(v)|$ for each $v \in S$, and the observation of the previous paragraph shows that the Hadamard product $h = f * g$ satisfies $h_v^{(d)}(h) < d|k(v)|$ for each $v \in S$.

We may now apply the criterion 3.2 to $h \in K[[t]]$ with $r := d$.  By the previous paragraph, the left-hand side of the condition~(\ref{criterion}) of the criterion is at least $\delta/4d^3$. By the assumption that $h$ is $A$-analytic, Lemma 3.3 shows that the right-hand side of~(\ref{criterion}) equals $(2d+1)\sum_v \log^+{R_v^{-1}}$. We thus obtain the algebraicity of $h$ as soon as $\log^+{R_v^{-1}} \leq \delta/12d^4 < \delta\big/(2d+1)4d^3$, giving us part (i) of Theorem 1.7.

For part (ii), we may apply the criterion 3.2 with $r = 1$. Under the stated additional assumptions the sequence $a(n)$ of coefficients of $f \in K[[t]] \cap K(t)$ has also pure period dividing $|k(v)^{\times}|$ at each $v \in S$, hence the Hadamard quotient in fact satisfies $h_v^{(1)}(h) < |k(v)|$. Thus 3.2 yields the rationality of the Hadamard quotient $h$ as soon as $h$ is $A$-analytic and $\sum_v \log^+{R_v^{-1}} < \delta/12$. \proofend

\bigskip

{\sc Appendix: A conjectural sharp form of the algebraicity criterion}

\bigskip

We propose a conjectural sharp form of the algebraicity criterion 3.2 used in the proof of Theorem 1.7. It generalizes at the same time a conjecture of I. Ruzsa~\cite{rusza,christol,zannierruzsa} and the  rationality criteria of Borel-Dwork, P\'olya-Bertrandias, and Bost-Chambert-Loir~\cite{bostacl}, which it unifies in a common framework. Partial results on this conjecture will appear in a subsequent paper.

\medskip

{\bf A.1.} With the notation and assumptions of 1.3, consider a $K$-rational point $P \in \proj^N(K)$ of a projective space and a smooth formal curve $\widehat{C}$ through $P$ in $\proj_K^N$, which is to say a one-dimensional smooth formal subscheme of the $P$-adic formal completion  of the scheme $\proj_K^N$. In the terminology of~\cite{bostacl} which inspires our Definition 1.6, assume $\widehat{C}$ is an {\it $A$-analytic curve}; we refer to 3.7 in {\it loc. cit.} for the general definition, but note that the condition is fulfilled if for a finite set of places $S$ the germ $\widehat{C}$ is $v$-adic analytic for each $v \in S$ and has a model over $O_{K,S}$. We also note that the graph of a formal function which is $A$-analytic in the sense of 1.4 is an $A$-analytic formal curve.

Consider Bost's adelic canonical semi-norm on the $K$-line $T_P\widehat{C}$, and its degree
$$
\widehat{\deg} \, (T_P\widehat{C},\|\cdot\|_v^{\mathrm{can}}) := \sum_v (-\log{\|\xi\|_v^{\mathrm{can}}})  \in \R \cup \{+\infty\}.
$$
Here, the summation is over all places $v$ of $K$ with the convention $-\log{0} := +\infty$; $\xi \in T_P\widehat{C} \, \setminus \, \{0\}$ is an arbitrary generator of the $K$-line $T_P\widehat{C}$; and the canonical $v$-adic semi-norms $\|\cdot\|_v^{\mathrm{can}}$ are defined as in~\cite{bostacl} (see Def.~4.7 in~{\it loc. cit.}). We refer to the last paper, as well as to Bost's articles~\cite{bostgerms} and~\cite{bosteval}, also for the close relationship of the $v$-adic canonical semi-norm with the analytic notion of capacity and the purely metric notion of transfinite diameter, both with respect to the divisor $P \in \widehat{C}$, of the domain in $\proj_{K_v}^N$ of $v$-adic analyticity of $\widehat{C}$. We have used this relationship implicitly in the proof of Theorem 1.4. The authors of~\cite{bostacl} only consider the number field case, but the situation is unaltered in the function field case.

It is proved in~\cite{bostgerms,bosteval,bostacl} that the positivity condition $\widehat{\deg} \, (T_P\widehat{C},\|\cdot\|_v^{\mathrm{can}}) > 0$ implies that $\widehat{C}$ is (the germ of) an algebraic subscheme of $\proj_K^N$. It is easy to see, on the other hand, that there are uncountably many possibilities for $\widehat{C}$ such that $\widehat{\deg} \, (T_P\widehat{C},\|\cdot\|_v^{\mathrm{can}}) = 0$.
 We propose the following as a common conjectural strengthening of the last result.
 For every closed point $s \in B$, let $h_s(\widehat{C}) \in \N_0 \cup \{+\infty\}$ be $+\infty$, if either $\widehat{C}$ does not have a reduction in $\proj_{k(s)}^N$ or if this reduction is not an algebraic subscheme; and the degree of the mod $s$ reduction $\widehat{C}_s$ as an algebraic curve of $\proj_{k(s)}^N$, otherwise. Our conjecture is then:

\medskip

\begin{quote}
{\it  If
$$
\widehat{\deg} \, (T_P\widehat{C},\|\cdot\|_v^{\mathrm{can}}) + \liminf_{n \in \N} \, \Big\{ \frac{1}{n} \sum_{s : \, h_s(\widehat{C}) < n} \log{|k(s)|} \Big\} > 0,
$$
 then the $A$-analytic curve $\widehat{C}$ is an algebraic subscheme of $\proj_K^N$.
}
\end{quote}

\medskip

{\bf A.2.} We  detail a particular case of Conjecture~A.1 as it applies to the rationality of a formal function on a regular, projective, and geometrically connected curve $X/K$. It implies in particular Ruzsa's conjecture from~\cite{rusza}, in the stronger form considered by Perelli-Zannier~\cite{perellizannier,zannierruzsa}  and Christol~\cite{christol}: the subring $\Z[1/(1-t)] \subset \Z[[t]]$ is characterized by the conditions that (1) the radius of convergence is strictly $> 1/e$; and (2) for a set of primes $p$ of full density, the mod $p$ reduction is of the form $A_p/(1-x)^p$ with $A_p \in \F_p[t]$ a polynomial of degree $< p$. (Our conjecture below in fact predicts a yet stronger form of Ruzsa's original conjecture, in which the condition on the radius of convergence is replaced by an assumption of meromorphy on a simply connected domain $G \ni 0$ such that the conformal mapping radius of $(G,0)$ strictly exceeds $1/e$.)

Let $S$ be a finite non-empty set of places of $K$, assumed to contain all the places of bad reduction of $X$, as well as (in the number field case)  all the archimedean places. As before, fix a $K$-rational point $P \in X(K)$, and consider a local parameter $t \in K(X)$ at $P$. For each place $v$ of $K$ consider a non-constant rational function $\varphi_v \in K_v(X_v)$ regular outside $P$, and let $\Omega_v \subset X_v^{\mathrm{an}}$ be the complement of the lemniscate domain $\{ x \in X_v^{\mathrm{an}} \mid |\varphi_v(x)|_v \leq 1 \}$. Assume furthermore that  for all $v \notin S$ the set $\Omega_v$ is ``the unit disc'' $|t(x)|_v < 1$. With $-m_v \in \Z^{< 0}$ the order of $\varphi_v$ at $P$, write $\varphi_v = c_v t^{-m_v} + \cdots$ around $P$, with $c_v \in K_v \setminus \{0\}$. Then the real quantity $K_{\Omega,P} := \sum_v \log{|c_v|_v}$ depends only on the adelic set $\Omega = (\Omega)_v$ and the point $P \in X(K)$, but not on the choice of local parameter $t$ at $P$. It equals the {\it capacitary degree} $\widehat{\deg}(T_PX,\|\cdot\|_{\Omega}^{\mathrm{cap}}) \geq \widehat{\deg}(T_PX,\|\cdot\|_{\Omega}^{\mathrm{can}})$ from~\cite{bostacl}; it also coincides with the negated logarithm of either of Rumely's or Chinburg's capacities, with respect to the divisor $P \in X(K)$, of the adelic subset of $X(\A_K)$ whose $v$-adic component is the complement in $X(K_v)$ of $\Omega_v$.  One may consider either of these notions for adelic sets more general than lemniscate domains, to which for simplicity we restrict here.

\medskip

{\bf Conjecture. }
{\it (i) Let $f \in \widehat{\Oz_{X,P}}$ be a formal function around $P$ defined over $O_{K,S}$. For each $s \in B$ let $d_f(s) \in \N_0 \cup \{+\infty\}$ be $+\infty$ if either $s \in S$ or the reduction $f \mod{s} \in \widehat{\Oz_{X_s,P_s}}$ is not in $k(s)(X_s)$; otherwise define it to be the degree of the reduction, viewed as a rational map $X_s \to \proj_{k(s)}^1$.

 If
$$
K_{\Omega,P} + \liminf_{n \in \N} \, \Big\{ \frac{1}{n} \sum_{s : \, d_f(s) < n} \log{|k(s)|} \Big\} > 0
$$
and $f$ extends for each $v \in S$ to a meromorphic function on $\Omega_v$, then $f \in K(X)$ is a rational function.

(ii) On the other hand, for any adelic set $\Omega$ as before, and any function $h : B \to \N_0 \cup \{+\infty\}$ such that
$$
K_{\Omega,P} + \liminf_{n \in \N} \, \Big\{ \frac{1}{n} \sum_{s : \, h(s) < n} \log{|k(s)|} \Big\} = 0,
$$
there are uncountably many formal functions $f \in \widehat{\Oz_{X,P}}$ around $P$ which are defined over $O_{K,S}$, analytic on $\Omega$, and satisfy $d_f(s) = h(s)$ for each $s \notin S$.}

\medskip

For example, if $X = \proj_K^1$, $P = 0$, and $\Omega_v$ is the disc $|x|_v < R_v$, where $R_v = 1$ for all $v \notin S$, then $K_{\Omega,P} = \sum_v \log{R_v}$. In this case part (ii) of the conjecture is obvious, while a simple sieving argument  shows that an $f \in O_{K,S}[[x]]$ satisfying the inequality of part~(i) is completely determined by a sufficiently long polynomial truncation; in particular, the power series in (i) form a countable set. This observation, together with the classical theorem of P\'olya and Bertrandias (see~\cite{amice}) and its generalization to higher genus in~\cite{bostacl} by Bost and Chambert-Loir, is our main motivation in support of the conjectural generalization stated here of the algebraicity criterion from~3.2.

\end{document}